\documentclass[conference]{IEEEtran}
\IEEEoverridecommandlockouts

\usepackage{enumitem}

\usepackage{hyperref}
\usepackage{subcaption}
\usepackage{graphicx}
\usepackage[version=4]{mhchem}
\usepackage{siunitx}
\usepackage{algorithm}
\usepackage{mathrsfs}
\usepackage{cite}
\usepackage{amsmath,amssymb,amsfonts,bm}
\usepackage{textcomp}
\usepackage{xcolor}
\usepackage{longtable,tabularx}
\usepackage[bottom= 54pt, right= 54pt, left= 54pt, top= 54pt]{geometry}
\usepackage[noend]{algpseudocode}
\begin{document}

\title{\vspace*{18pt} Investigating the Performance of EKF, UKF, and PF for Quadrotor Position Estimation in Hurricane Wind Disturbances}



\author{Ahmed A. Elgohary, Benjamin Gwinnell, and Josh Augustine} 



\maketitle
\begin{abstract} 
Natural disasters, such as hurricanes and typhoons, pose significant challenges to public safety and infrastructure. While government agencies rely on multi-million-dollar UAV systems for storm data collection and disaster response, smaller drones lack the ability to autonomously adapt to rapidly changing environmental conditions, such as turbulent winds. This project investigates the implementation of advanced state estimation filters—Extended Kalman Filter (EKF), Unscented Kalman Filter (UKF), and Particle Filter (PF)—to enhance the control and adaptability of quadrotor drones in uncertain wind profiles modeled by Von Karman turbulence. While the EKF relies on linearization techniques using the Taylor series expansion, which can struggle under high nonlinearity, the UKF leverages sigma points for better performance in nonlinear systems without requiring Jacobian computations. The PF, on the other hand, addresses non-Gaussian noise and severe nonlinearities by employing a large number of particles, albeit at a high computational cost. To enhance accuracy and minimize estimation errors, a genetic algorithm (GA) was employed to optimally tune the process and measurement noise covariances (\(Q\) and \(R\) matrices) as well as UKF-specific parameters (\(\alpha\), \(\beta\), \(\kappa\)). This study highlights the trade-off between accuracy, computational efficiency, and smoothing capabilities across these filters. Despite its robustness, the PF suffered from computational inefficiencies due to the high state dimensionality, while the EKF demonstrated faster computation but lower adaptability in nonlinear conditions. The UKF emerged as a balanced approach, achieving superior performance in capturing dynamic wind disturbances.

\end{abstract}

\section{Introduction}
Unmanned aerial vehicles (UAVs), including quadrotors, have become indispensable tools in various fields, ranging from environmental monitoring and disaster relief to search and rescue operations \cite{NOAA2020}. Despite their increasing affordability and versatility, publicly available drones are often ill-equipped to handle extreme environmental conditions, such as hurricane-force winds, compared to the multi-million-dollar UAVs used by organizations like NOAA. These government-grade systems, while highly effective, are costly and not widely accessible. The challenge lies in enabling smaller, commercially available drones to perform similarly under adverse conditions, ensuring reliable operation and accurate positional control despite turbulent wind profiles. Turbulence models such as Von Karman has been widely adopted to simulate hurricane-level disturbances, enabling the evaluation of control system robustness \cite{VonKarman1948}. State estimation filters, such as EKF, UKF, and PF, have proven effective in handling noisy and uncertain environments, building on Kalman’s seminal work in optimal filtering \cite{Kalman1960}. This project seeks to bridge the gap between commercially available drones and government-grade UAVs by incorporating advanced state estimation filters into quadrotor control systems. Utilizing turbulence models such as Von Karman, and optimizing filter parameters, including process and measurement noise covariances, via Genetic Algorithms \cite{Goldberg1989}, this study presents a thorough evaluation of filter performance under hurricane-level wind conditions. The findings offer actionable insights into enhancing the reliability, accuracy, and safety of quadrotor operations in dynamic and challenging environments.

\section{Background and Preliminaries}
In this section, we present the background and preliminaries relevant to this study, including the wind disturbance model used in the simulation and the quadrotor dynamics model and equations. Additionally, we briefly talk about the background of the algorithms used during the simulation results.

\subsection{Wind Disturbance Model}
To accurately simulate hurricane-force disturbances and evaluate the implementation of state estimation filters in quadrotor control systems, wind turbulence models play a pivotal role. This study employs a widely recognized model: the Von Karman Wind Turbulence Model. This model is extensively utilized by institutions such as the U.S. Department of Defense and NASA for simulating wind disturbances in aircraft systems \cite{nasa}. The mathematical representations of this model provide realistic testbeds for assessing the robustness of control strategies. The Von Karman wind turbulence model introduces second-order differential equations to represent the longitudinal, lateral, and vertical wind components:
\begin{align}
    \frac{d^2 u_g}{dt^2} + 2 \zeta_u \omega_{n_u} \frac{du_g}{dt} + \omega_{n_u}^2 u_g &= \omega_{n_u}^2 \sigma_u^2 w_u, \\
    \frac{d^2 v_g}{dt^2} + 2 \zeta_v \omega_{n_v} \frac{dv_g}{dt} + \omega_{n_v}^2 v_g &= \omega_{n_v}^2 \sigma_v^2 w_v, \\
    \frac{d^2 w_g}{dt^2} + 2 \zeta_w \omega_{n_w} \frac{dw_g}{dt} + \omega_{n_w}^2 w_g &= \omega_{n_w}^2 \sigma_w^2 w_w.
\end{align} 

Where \(u_g\), \(v_g\), and \(w_g\) denote the longitudinal, lateral, and vertical wind velocity components, respectively. The parameters \(\zeta_u\), \(\zeta_v\), and \(\zeta_w\) represent the damping ratios for the longitudinal, lateral, and vertical wind components, while \(\omega_{n_u}\), \(\omega_{n_v}\), and \(\omega_{n_w}\) are the natural frequencies of turbulence in their respective directions. The variables \(\sigma_u\), \(\sigma_v\), and \(\sigma_w\) correspond to the turbulence intensities (standard deviations) in each direction, and \(w_u\), \(w_v\), and \(w_w\) are the white noise inputs driving the turbulence model. In this study, only the longitudinal effects are taken into consideration, specifically the effects of \(u_g\) and \(w_g\).

\subsection{Quadrotor Dynamics}
Quadrotor dynamics are inherently nonlinear and are typically modeled in six degrees of freedom, encompassing translational and rotational motions. These dynamics are influenced by aerodynamic forces, gravitational forces, and control inputs, such as thrust and torque. To simplify the analysis and focus on the vertical and forward motion of the quadrotor, this study employs the longitudinal dynamics model as presented in \cite{uav_book}. The longitudinal dynamics of the quadrotor are represented in a state-space form as:

\begin{equation}
    \dot{\mathbf{x}} = f(\mathbf{x}) + g(\mathbf{x}) \mathbf{u},
\end{equation}

where the state vector \(\mathbf{x}\), system dynamics \(f(\mathbf{x})\), control input matrix \(g(\mathbf{x})\), and control input vector \(\mathbf{u}\) are defined as follows:

\[
\mathbf{x} =
\begin{bmatrix}
    p_n \\
    h \\
    u \\
    w \\
    \theta \\
    q
\end{bmatrix},
\quad
f(\mathbf{x}) =
\begin{bmatrix}
    u \cos \theta + w \sin \theta \\
    u \sin \theta - w \cos \theta \\
    -qw - g \sin \theta \\
    qu + g \cos \theta \\
    q \\
    0
\end{bmatrix},
\]

\[
g(\mathbf{x}) =
\begin{bmatrix}
    0 & 0 \\
    0 & 0 \\
    -\frac{1}{m} & 0 \\
    0 & 0 \\
    0 & 0 \\
    0 & \frac{1}{J_y}
\end{bmatrix},
\quad
\mathbf{u} =
\begin{bmatrix}
    F \\
    \tau_\theta
\end{bmatrix}.
\]

The state vector \(\mathbf{x}\) consists of the northward position \(p_n\), altitude \(h\), forward velocity \(u\), vertical velocity \(w\), pitch angle \(\theta\), and pitch rate \(q\). The control input vector \(\mathbf{u}\) includes \(F\), the thrust force, and \(\tau_\theta\), the torque about the pitch axis. The parameters \(m\), \(g\), and \(J_y\) represent the mass of the quadrotor, gravitational acceleration, and pitch-axis moment of inertia, respectively. These equations provide a simplified yet accurate representation of the quadrotor's motion, enabling effective implementation and testing of control and state estimation strategies.

\subsection{State Estimation and Optimization Algorithms}
The EKF extends the classical Kalman Filter to handle nonlinear systems by linearizing the dynamics and measurement models around the current estimate using a Taylor series expansion. While computationally efficient, the EKF is sensitive to high levels of nonlinearity and relies heavily on accurate modeling of the system dynamics. The EKF algorithm used in this study is implemented based on the approach described in \cite{EKF}. The UKF improves upon the EKF by eliminating the need for Jacobian matrix computations. Instead, the UKF leverages sigma points to capture the mean and covariance of the state distribution more accurately, especially in highly nonlinear systems. This makes the UKF particularly effective for handling nonlinear dynamics without the limitations of linearization. The implementation of the UKF in this study follows the methodology outlined in \cite{UKF}. The PF is a powerful nonparametric filter that uses a set of weighted particles to approximate the posterior distribution of the state. Unlike the EKF and UKF, the PF can handle severe nonlinearity and non-Gaussian noise. However, its computational cost increases significantly with the dimensionality of the state space, requiring a large number of particles for accurate estimation. The algorithm for the PF used in this study is derived from \cite{PF}. Additionally, the Genetic Algorithm (GA) is employed for optimizing the process and measurement noise covariance matrices. The MATLAB and Simulink Optimization Toolbox is utilized for implementing the GA, which plays a critical role in minimizing the estimation error and ensuring a balance between accuracy and computational efficiency during the simulation results \cite{matlab_ga}.

\section{Main Results}
In this section, we present the parameters used for simulating the quadrotor dynamics and wind turbulence models, as well as the results obtained for hovering and wind disturbance conditions. Additionally, we will present the comparative results for the position estimation of the quadrotor under wind disturbance using the three filters (EKF, UKF, PF) and the Genetic Algorithm (GA) with different cost functions.

\subsection{Quadrotor Model in Hovering Condition}

\begin{figure}[ht]
    \centering
    \includegraphics[width=1\linewidth]{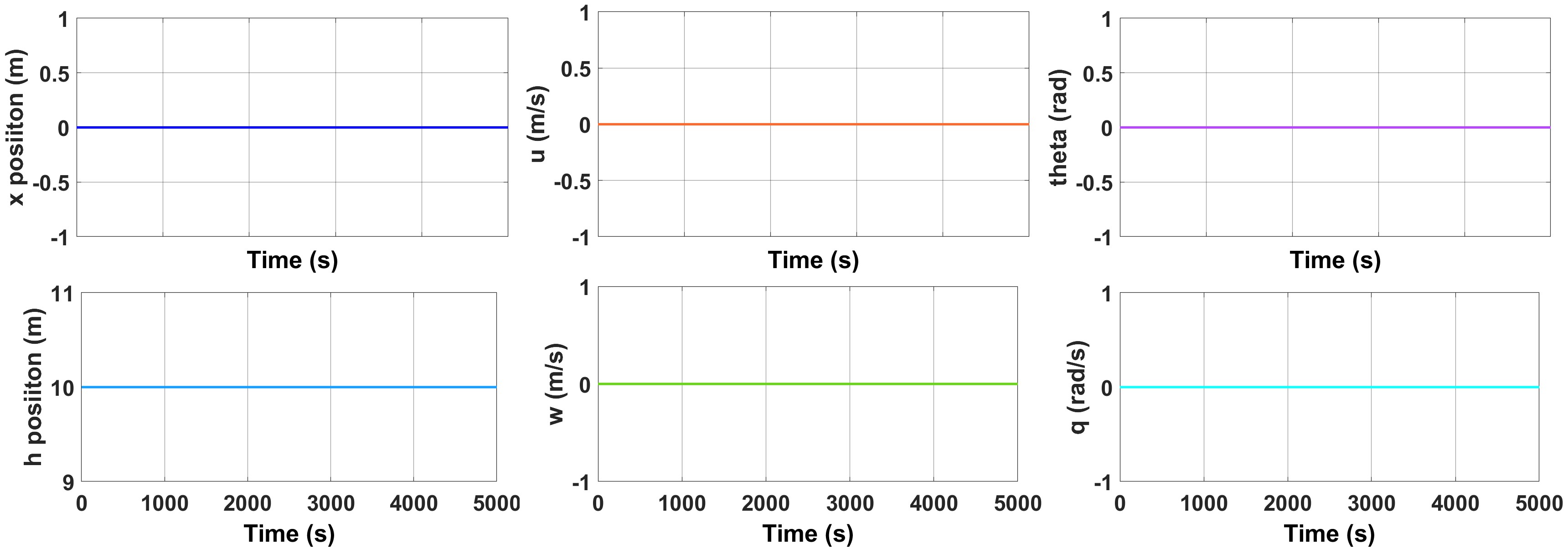}
    \caption{Quadrotor hovering condition at 10 m altitude with steady-state variables}
    \label{fig:hovering}
\end{figure}

The quadrotor was simulated during hovering at an altitude of 10 m using the following parameters: mass \(m = 1.5 \, \mathrm{kg}\), gravitational acceleration \(g = 9.81 \, \mathrm{m/s^2}\), moment of inertia about the pitch axis \(J_y = 0.057 \, \mathrm{kg \cdot m^2}\), control inputs \(\tau_\theta = 0\), and \(F = m g\) to maintain hovering. After applying these parameters, the quadrotor model was simulated for 5000 seconds. The simulation results demonstrated that all state variables remained constant over time, indicating that the quadrotor successfully maintained the hovering condition. This is illustrated in Figure \ref{fig:hovering}, where the states remain steady throughout the simulation, validating the hovering condition.

\subsection{Wind Turbulence Model and its Effect on Quadrotor Dynamics}

\begin{figure}[ht]
    \centering
    \includegraphics[width=1\linewidth]{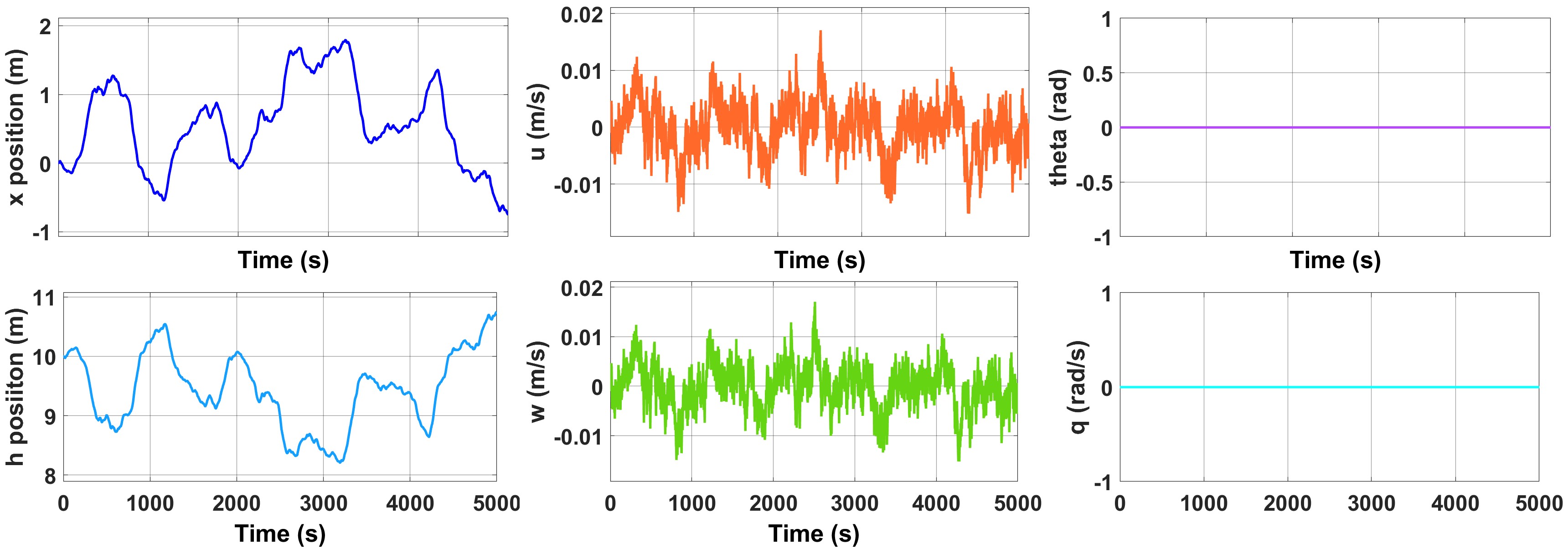}
    \caption{Quadrotor response to wind disturbances showing deviations in altitude and position from the hovering condition}
    \label{fig:hover_dis}
\end{figure}

The Von Karman wind turbulence model was used to introduce wind disturbances into the simulation. The following parameters were utilized for the wind turbulence model: altitude \(10 \, \mathrm{m}\) and wind speed \(70 \, \mathrm{m/s}\) (corresponding to a Category 5 hurricane). The longitudinal and vertical wind disturbances, \(u_g\) and \(w_g\), were added to the quadrotor dynamics during hovering. Under these conditions, the quadrotor states responded dynamically, as shown in Figure \ref{fig:hover_dis}. The altitude (\(h\)) and northward position (\(p_n\)) exhibited significant deviations from the steady-state hovering condition due to the wind disturbances, demonstrating the impact of the turbulence on the quadrotor dynamics.

\subsection{Quadrotor Dynamics and Filter Setup}

In this subsection, we construct the quadrotor dynamics model and prepare it for implementation across all filters. Additionally, we present the comparative filter results obtained using the optimal parameters identified through the Genetic Algorithm (GA). The state vector is defined as:
\[
x_k =
\begin{bmatrix}
    p_n \\
    h \\
    u \\
    w \\
    \theta \\
    q
\end{bmatrix},
\]
where \(p_n\) is the northward position, \(h\) is the altitude, \(u\) is the forward velocity, \(w\) is the vertical velocity, \(\theta\) is the pitch angle, and \(q\) is the pitch rate. The measurement vector is given as:
\[
z_k =
\begin{bmatrix}
    p_n \\
    h
\end{bmatrix},
\]
which represents the quadrotor's position under the effects of wind disturbance and Gaussian noise. The noise is added to the position measurements with zero mean and a covariance of 0.01 throughout the entire simulation duration.

For the Extended Kalman Filter (EKF), as outlined in \cite{EKF}, the algorithm relies on linearizing the nonlinear dynamics of the system. During the error covariance estimation step, which is the second step in the EKF algorithm, the Jacobian matrix \(F_k\) of the system dynamics is computed as:
{\small
\[
F_k =
\frac{\partial f(x, u)}{\partial x} =\\
\begin{bmatrix}
    0 & 0 & \cos\theta & \sin\theta & -u\sin\theta + w\cos\theta & 0 \\
    0 & 0 & \sin\theta & -\cos\theta & u\cos\theta + w\sin\theta & 0 \\
    0 & 0 & 0 & -q & -g\cos\theta & -w \\
    0 & 0 & q & 0 & -g\sin\theta & u \\
    0 & 0 & 0 & 0 & 0 & 1 \\
    0 & 0 & 0 & 0 & 0 & 0
\end{bmatrix}
\]
}
Similarly, the measurement model requires the Jacobian matrix \(H_k\), defined as:
\[
H_k =
\frac{\partial h(x)}{\partial x} =
\begin{bmatrix}
    1 & 0 & 0 & 0 & 0 & 0 \\
    0 & 1 & 0 & 0 & 0 & 0
\end{bmatrix}.
\]

These matrices, along with the state dynamics and measurement noise, are integral to the implementation of the EKF algorithm. The setup forms the basis for comparing the performance of all filters under varying wind disturbance conditions. The process noise covariance \(Q_k\) and measurement noise covariance \(R_k\) used in the simulations are defined as diagonal matrices, as shown below:
\[
Q_k =
\begin{bmatrix}
    q_1 & 0 & 0 & 0 & 0 & 0 \\
    0 & q_2 & 0 & 0 & 0 & 0 \\
    0 & 0 & q_3 & 0 & 0 & 0 \\
    0 & 0 & 0 & q_4 & 0 & 0 \\
    0 & 0 & 0 & 0 & q_5 & 0 \\
    0 & 0 & 0 & 0 & 0 & q_6
\end{bmatrix}, \quad
R_k =
\begin{bmatrix}
    r_1 & 0 \\
    0 & r_2
\end{bmatrix}.
\]

To optimize the values of the parameters in \(Q_k\) and \(R_k\), the Genetic Algorithm (GA) provided by MATLAB's Simulink toolbox was employed. The optimization was based on minimizing the error between the true position and the estimated position obtained from the EKF. The cost function used for this optimization is the Least Squares Method, represented as:
\[
J = \sum_{i=1}^n (x_{\text{true},i} - \hat{x}_i)^2 + \sum_{i=1}^n (h_{\text{true},i} - \hat{h}_i)^2.
\]

\begin{figure}[ht]
    \centering
    \includegraphics[width=1\linewidth]{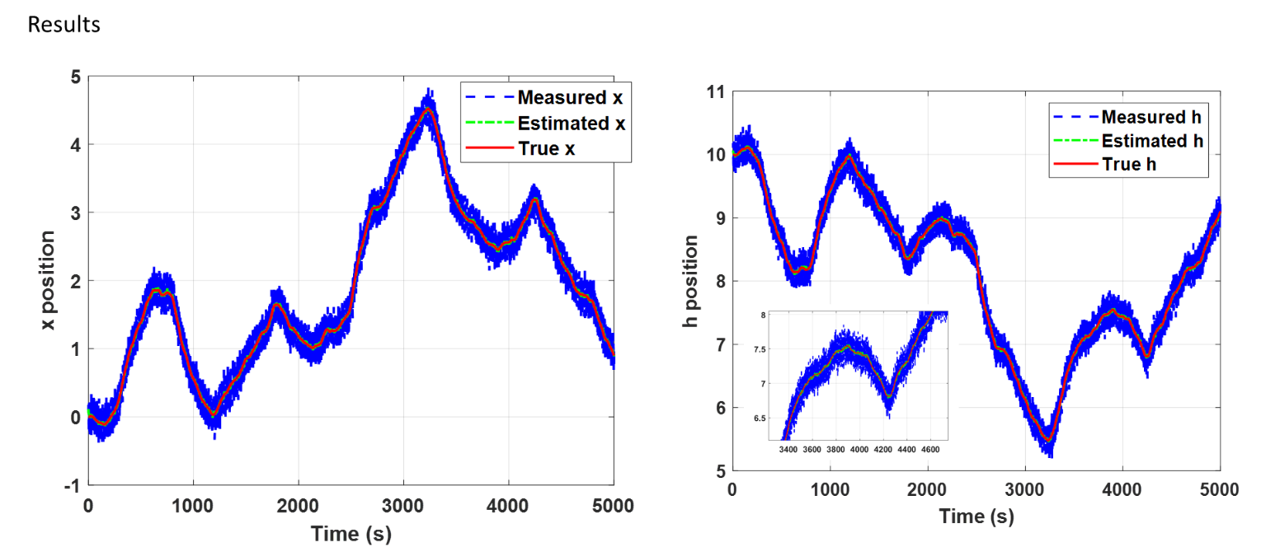}
    \caption{Comparison of true, estimated by EKF, and measured positions}
    \label{fig:EKF}
\end{figure}

Where \(x_{\text{true},i}\) and \(\hat{x}_i\) represent the true and estimated values of the northward position at the \(i\)-th time step, respectively. Similarly, \(h_{\text{true},i}\) and \(\hat{h}_i\) correspond to the true and estimated values of the altitude at the \(i\)-th time step. The GA parameters used throughout the simulations, including for the UKF, were as follows: maximum number of generations = 600, population size = 50, mutation rate = 0.167, crossover probability = 0.8, and stop condition = \(1e{-6}\). Using the GA as mentioned earlier, we obtained the optimal values of \( Q \) and \( R \) as follows: \( R = [0.9999, 0.9999] \) and \( Q = [4.83e\text{-}05, 2.25e\text{-}03, 1.25e\text{-}03, 1.86e\text{-}03, 2.23e\text{-}05, 2.19e\text{-}05]. \) Using these optimal values, we generated comparative results for the true values, estimated values, and measured values for the positions \( x \) (representing \( p_n \)) and \( h \). The results are illustrated in Figure \ref{fig:EKF}.

In addition, we sought to further demonstrate the effectiveness of the Genetic Algorithm (GA) by adjusting the smoothing of the position estimation in cases where the true data is unavailable for calculating \(J\). To achieve this, a different cost function composed of four terms with corresponding weights was utilized. The first two terms account for the accuracy of the estimation between the measured and estimated data, while the last two terms emphasize the smoothing of the estimated positions. The cost function \(J\) is defined as follows:

{\small
\[
J = \beta_1 \cdot \sqrt{\frac{1}{n} \sum_{i=1}^n (x_{\text{measured},i} - \hat{x}_i)^2} + 
\beta_2 \cdot \sqrt{\frac{1}{n} \sum_{i=1}^n (h_{\text{measured},i} - \hat{h}_i)^2} + 
\]
\[
\alpha_1 \cdot \sum |\text{diff}(\hat{x})| + 
\alpha_2 \cdot \sum |\text{diff}(\hat{h})|.
\]
}

\begin{figure}[ht]
    \centering
    \includegraphics[width=1\linewidth]{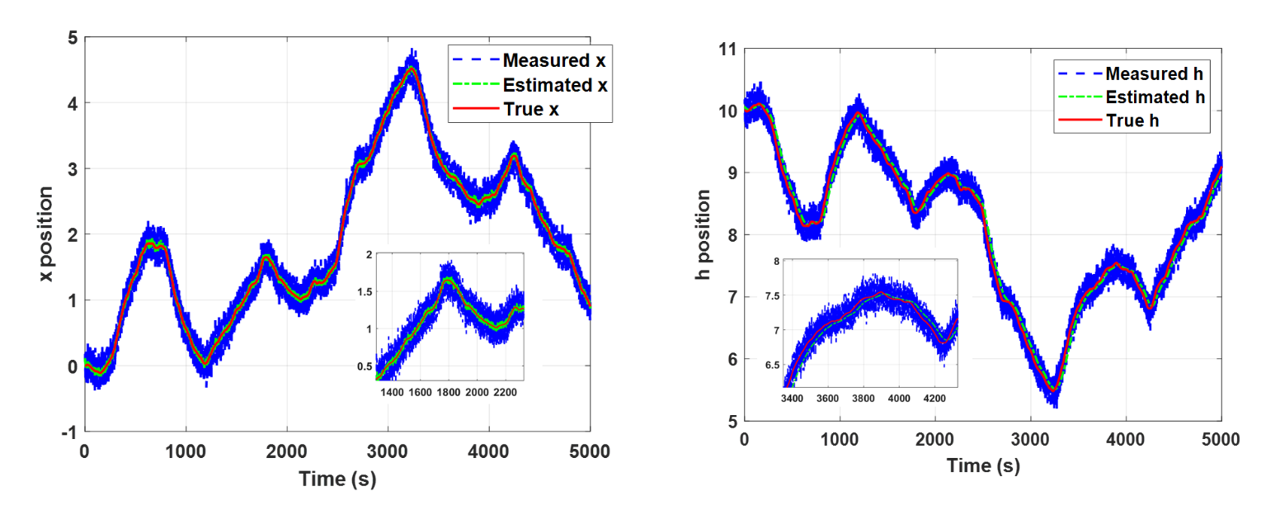}
    \caption{Comparison of true and estimated positions by EKF using the optimized weights \(\beta_1 = 1\), \(\beta_2 = 60\), \(\alpha_1 = 10\), \(\alpha_2 = 5\), highlighting smoother and more accurate \(z\) and \(x\) estimations}
    \label{fig:EKF_smooth}
\end{figure}

Where the weights used in this case are \(\beta_1 = 1\), \(\beta_2 = 60\), \(\alpha_1 = 10\), and \(\alpha_2 = 5\), balancing the contributions of each term in the cost function. These values were chosen to achieve an optimal trade-off between estimation accuracy and smoothness. Using this modified cost function, the GA optimized the noise covariance matrices as \(R = \text{diag}(0.9999, 0.999)\) and \(Q = \text{diag}(9.9924\text{e-05}, 3.7216\text{e-05}, 7.5426\text{e-05}, 2.9733\text{e-04},\\
9.3899\text{e-05}, 1.7283\text{e-05})\). Using these optimal \(Q\) and \(R\) matrices, we obtained the following results for the EKF in position estimation. The true values and estimated values by EKF are shown in Figure \ref{fig:EKF_smooth}, where, as highlighted in the zoomed-in section, the \(z\) and \(x\) estimation is smoother and closer to the true value. This demonstrates and validates our approach. This framework proved effective in achieving accurate and smooth position estimations, even in the absence of true data.

\begin{figure}[ht]
    \centering
    \includegraphics[width=1\linewidth]{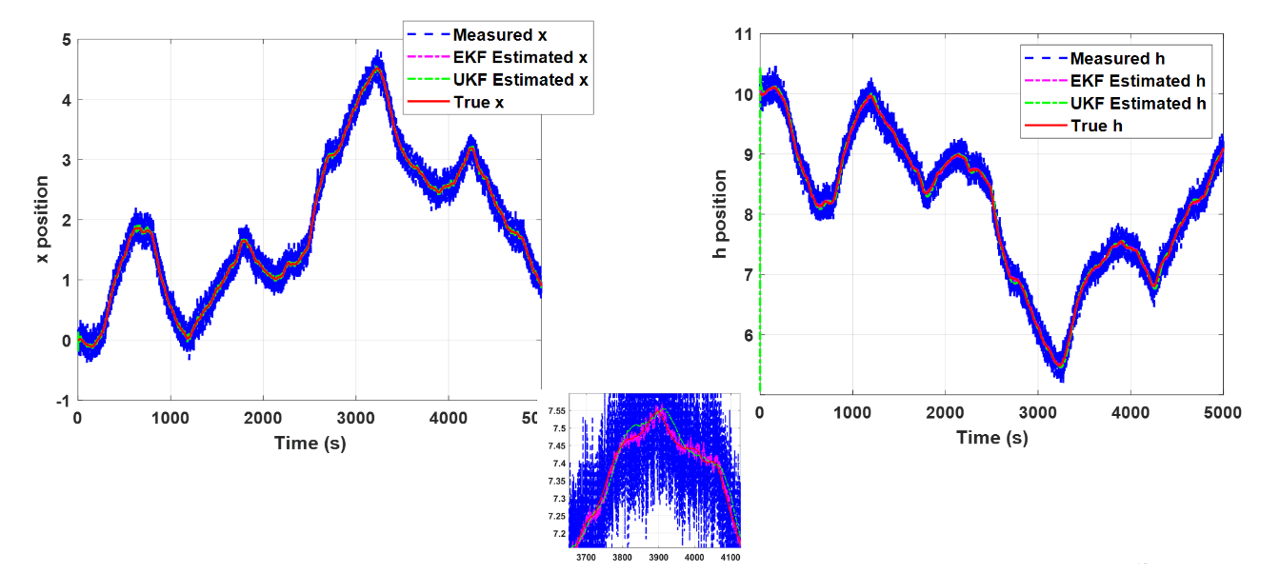}
    \caption{Comparison of position estimations by EKF and UKF, highlighting the UKF's smoother and more accurate performance in the zoomed section}
    \label{fig:UKF}
\end{figure}

Using the UKF algorithm based on \cite{UKF}, which assigns sigma points to the nonlinear function without requiring the Jacobian matrix, three additional parameters—\(\alpha\), \(\beta\), and \(\kappa\)—are introduced. Using the Genetic Algorithm (GA) to estimate these parameters alongside the \(Q\) and \(R\) matrices, and employing the same cost function as the EKF, which minimizes the error between the true and estimated positions, we obtained the following optimal parameters: \(\alpha = 0.4\), \(\beta = 2\), \(\kappa = 0\), \(R = \text{diag}(0.9999999, 0.9999999)\), and \(Q = \text{diag}(1.0\text{e-18}, 1.0\text{e-11}, 1.0\text{e-19}, 1.0\text{e-10}, 1.0\text{e-17}, 1.0\text{e-12})\). Using these optimal values, we plotted in Figure~\ref{fig:UKF} the comparison between the position estimation for the measured, true, and estimated values using both EKF and UKF. As shown in the zoomed part, the UKF demonstrated a more precise and smoother estimation for the position compared to the EKF estimation.

\begin{figure}[ht]
    \centering
    \includegraphics[width=1\linewidth]{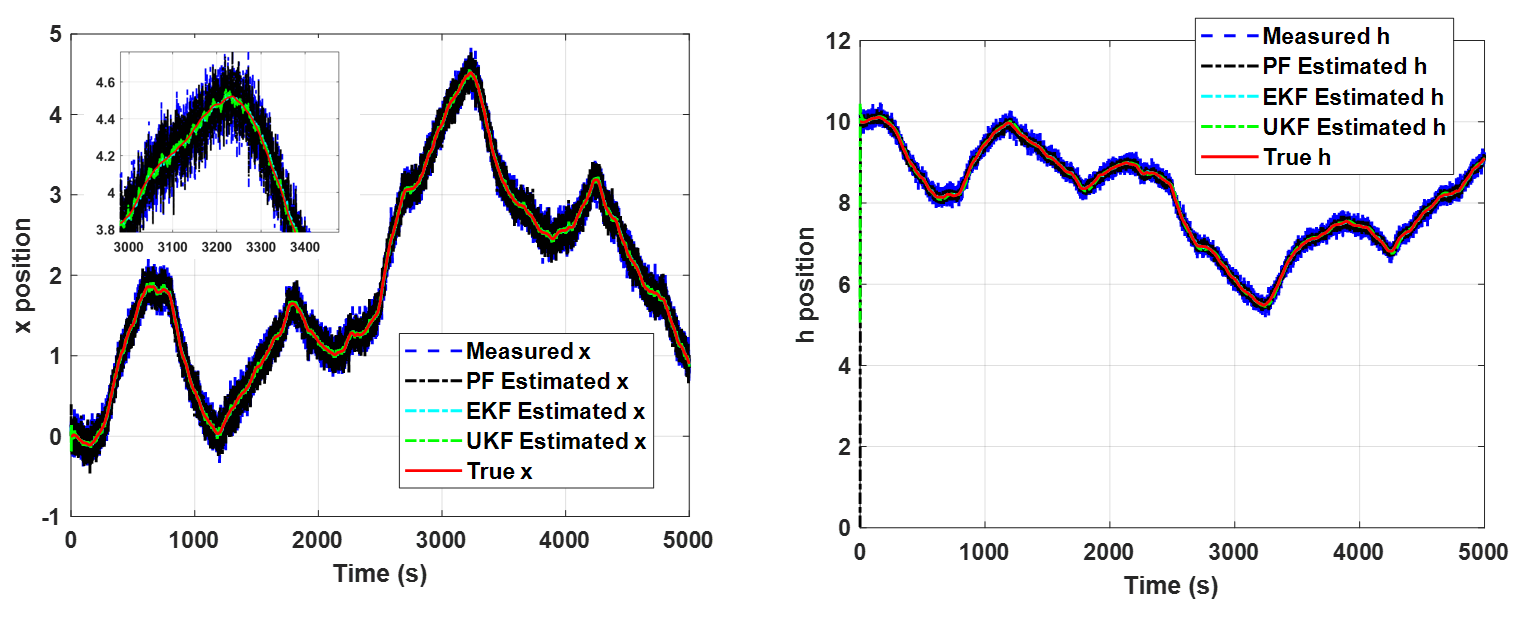}
    \caption{PF}
    \label{fig:PF}
\end{figure}

We also implemented the Particle Filter (PF) based on \cite{PF}, where a large number of particles (\(N\)) were assigned to approximate the posterior distribution of the states. The tuning parameters for the PF primarily include \(N\), which represents the number of particles, as well as the process (\(Q\)) and measurement (\(R\)) noise covariance matrices. Due to the high computational demand of the PF and the slow simulation time required for our quadrotor model with six states, we could not tune the PF parameters using the GA. Instead, we manually tuned these parameters, which resulted in suboptimal performance. Using 5000 particles and the following covariance matrices, \(R = \text{diag}(0.9999, 0.999)\) and \(Q = \text{diag}(4.8254\text{e-05}, 2.2522\text{e-03}, 1.2468\text{e-03}, 1.8631\text{e-03},\\ 2.2304\text{e-05}, 2.1894\text{e-05})\), we compared the results between the EKF, UKF, and PF for position estimation. As shown in Figure~\ref{fig:PF}, the PF failed to provide a smooth and accurate estimation of the position due to the limitations of manual tuning, highlighting the challenges of applying PF for high-dimensional systems.

\section{Conclusive Remarks and Future Work}

In this study, we addressed the challenges of position estimation and state prediction for quadrotors operating under extreme wind disturbances, such as hurricane-level turbulence. By employing advanced state estimation algorithms, namely the Extended Kalman Filter (EKF), Unscented Kalman Filter (UKF), and Particle Filter (PF), we demonstrated the effectiveness of these techniques in handling the nonlinear dynamics and stochastic nature of quadrotor systems. Using Genetic Algorithms (GA), we optimized the process and measurement noise covariance matrices (\(Q\) and \(R\)), along with UKF-specific parameters (\(\alpha\), \(\beta\), and \(\kappa\)), to achieve enhanced estimation accuracy and smoothness. The UKF emerged as the most robust filter, achieving superior performance in capturing the dynamic behavior of the quadrotor, particularly in highly nonlinear scenarios. However, the PF, despite its ability to handle non-Gaussian noise and severe nonlinearities, faced computational challenges and underperformed due to manual parameter tuning. The results showed that state estimation filters, when optimized, can significantly improve quadrotor performance and reliability under challenging conditions. These findings underscore the importance of integrating advanced optimization techniques, such as GA, with estimation algorithms to enhance the robustness of UAV systems in real-world applications.

\subsection*{Future Work}

Future work will focus on overcoming the computational limitations of the Particle Filter by exploring adaptive resampling techniques and parallel processing methods. Additionally, extending the quadrotor dynamics model to include full six-degree-of-freedom motion and incorporating additional sensors, such as inertial measurement units (IMUs) and barometers, will enable more accurate and comprehensive state estimation. Another promising avenue is to investigate hybrid filtering approaches, combining the strengths of the UKF and PF, to achieve better trade-offs between computational efficiency and estimation accuracy. Furthermore, real-world flight experiments will be conducted to validate the proposed methodologies and assess their practical feasibility in disaster relief and environmental monitoring scenarios. Finally, incorporating more complex turbulence models, such as those with time-varying parameters, will allow for more realistic testing of quadrotor performance in dynamic and uncertain environments.

\bibliography{references}

\begin{thebibliography}{10}
\providecommand{\url}[1]{#1}
\csname url@samestyle\endcsname
\providecommand{\newblock}{\relax}
\providecommand{\bibinfo}[2]{#2}
\providecommand{\BIBentrySTDinterwordspacing}{\spaceskip=0pt\relax}
\providecommand{\BIBentryALTinterwordstretchfactor}{4}
\providecommand{\BIBentryALTinterwordspacing}{\spaceskip=\fontdimen2\font plus
\BIBentryALTinterwordstretchfactor\fontdimen3\font minus \fontdimen4\font\relax}
\providecommand{\BIBforeignlanguage}[2]{{%
\expandafter\ifx\csname l@#1\endcsname\relax
\typeout{** WARNING: IEEEtran.bst: No hyphenation pattern has been}%
\typeout{** loaded for the language `#1'. Using the pattern for}%
\typeout{** the default language instead.}%
\else
\language=\csname l@#1\endcsname
\fi
#2}}
\providecommand{\BIBdecl}{\relax}
\BIBdecl

\bibitem{NOAA2020}
{NOAA}, ``Noaa’s role in environmental monitoring and disaster relief,'' \emph{National Oceanic and Atmospheric Administration Publications}, 2020, accessed online at \url{https://www.noaa.gov/}.

\bibitem{VonKarman1948}
T.~Von~Karman, ``Progress in the statistical theory of turbulence,'' \emph{Journal of Applied Physics}, vol.~19, pp. 1232--1244, 1948.

\bibitem{Kalman1960}
R.~E. Kalman, ``A new approach to linear filtering and prediction problems,'' \emph{Transactions of the ASME—Journal of Basic Engineering}, vol.~82, pp. 35--45, 1960.

\bibitem{Goldberg1989}
D.~E. Goldberg, \emph{Genetic Algorithms in Search, Optimization, and Machine Learning}.\hskip 1em plus 0.5em minus 0.4em\relax Addison-Wesley, 1989.

\bibitem{nasa}
F.~B. Leahy, ``Discrete gust model for launch vehicle assessments,'' in \emph{13th Conference on Aviation, Range and Aerospace Meteorology/American Meteorology Society}, 2008.

\bibitem{uav_book}
R.~W. Beard and T.~W. McLain, \emph{Small unmanned aircraft: Theory and practice}.\hskip 1em plus 0.5em minus 0.4em\relax Princeton university press, 2012.

\bibitem{EKF}
M.~I. Ribeiro, ``Kalman and extended kalman filters: Concept, derivation and properties,'' \emph{Institute for Systems and Robotics}, vol.~43, no.~46, pp. 3736--3741, 2004.

\bibitem{UKF}
E.~A. Wan and R.~Van Der~Merwe, ``The unscented kalman filter for nonlinear estimation,'' in \emph{Proceedings of the IEEE 2000 adaptive systems for signal processing, communications, and control symposium (Cat. No. 00EX373)}.\hskip 1em plus 0.5em minus 0.4em\relax Ieee, 2000, pp. 153--158.

\bibitem{PF}
J.~Elfring, E.~Torta, and R.~Van De~Molengraft, ``Particle filters: A hands-on tutorial,'' \emph{Sensors}, vol.~21, no.~2, p. 438, 2021.

\bibitem{matlab_ga}
\BIBentryALTinterwordspacing
{MathWorks}, ``Genetic algorithm - matlab \& simulink,'' 2024, accessed: 2024-12-12. [Online]. Available: \url{https://www.mathworks.com/discovery/genetic-algorithm.html}
\BIBentrySTDinterwordspacing

\end{thebibliography}
\bibliographystyle{IEEEtran}
\end{document}